\documentclass[american,12pt]{article}
\usepackage{amssymb,amsfonts,bm}
\usepackage{babel,amscd,xypic}
\usepackage{amsmath}
\usepackage{times}
\newtheorem{theorem}{Theorem}
\newtheorem{proposition}{Proposition}
\newtheorem{lemma}{Lemma}
\newtheorem{defn}{Definition}

\newtheorem{remark}{Remark}
\newcommand{\ep}{\hfill$\blacksquare$}

\begin{document}

\title{Noncommmutative Gelfand Duality for not necessarily unital $C^*$-algebras, Jordan Canonical form, and the existence of invariant subspaces}
\author{Mukul S. Patel\\ $<$patel@math.uga.edu$>$\\ Department of Mathematics, University of Georgia\\ Athens, GA  30602, USA}
\date{August 12, 2005}
\maketitle

\begin{abstract}
Gelfand-Naimark duality (Commutative $C^*$-algebras $\equiv$ Locally compact Hausdorff spaces) is extended to
\begin{center}$C^*$-algebras $\equiv$  Quotient maps on locally compact Hausdorff spaces.\end{center}
Using this duality, we give for an \emph{arbitrary} bounded operator on a complex Hilbert space of several dimensions, a functional calculus and the existence theorem for nontrivial invariant subspace.
\end{abstract}

\setcounter{secnumdepth}{4}
\setcounter{tocdepth}{1}

\section{\normalsize{INTRODUCTION}}

Connes attaches $C^*$-algebras to various quotient spaces arising in geometry \cite{connes}. Conversely, we assign a natural quotient map to any given $C^*$-algebra. Of course, for commutative algebras,
the Gelfand-Naimark theorem does the job:
\begin{theorem}[Gelfand-Naimark]\label{gelfand}
A commutative $C^*$-algebra  $A$ is naturally isomorphic to $C(P(A))$, the algebra of complex valued continuous functions vanishing at infinity on $P(A),$ the space of pure states of $A$.
\end{theorem}

Note that $P(A)$ is a locally compact Hausdorff space, and the theorem sets up a functorial equivalence between the category of commutative $C^*$-algebras on one hand and the category of locally compact Hausdorff spaces on the other hand. Also, this functor maps the subcategory consisting of compact Hausdorff spaces to the subcategory consisting of \emph{unital} commutative $C^*$-algebras.

Since the original theorem \cite{gelfand}, there have been several noncommutative generalizations in various directions \cite{1,2,5,4,7} with varying degree of success. Our generalization (Theorem \ref{main})  is implemented by identifying the \emph{natural} noncommutative analog of locally compact Hausdorff space---an equivalence relation on, or equivalently a quotient of, a locally compact Hausdorff space. This embraces most commonly occurring geometric situations on one hand, and all $C^*$-algebras on the other. The key to this quotient is the following trivial observation: \textit{A $C^*$-algebra is commutative if and only if
all its irreducible Gelfand-Naimark-Segal representations
are pair-wise inequivalent.} Thus, the noncommutativity of an algebra is completely captured
by the equivalence relation given by equivalence of irreducible GNS representations.

There have been studies of $C^*$-algebras via continuous functions on groupoids (See \cite{renault}, for example). The latter include equivalence relations as a special case. However, a larger algebra is needed to capture the whole situation. Our main result in this direction, Theorem \ref{main}, asserts that the algebra $A$ is canonically isomorphic to a certain algebra of regular Borel \emph{measures} on an equivalence relation $R(A)$. We then take this equivalence relation, or equivalently the quotient map it entails, as quantum space:

\begin{defn}[Quantum spaces]\label{quantum space} 
By  \textbf{quantum space (resp. compact quantum space)} we shall mean a quotient map \[q: X\twoheadrightarrow Y\] where $X$ is a locally compact (resp. compact) Hausdorff space.  Then a \textbf{quantum group (resp. semigroup, groupoid, etc.) space} will be a group object (resp. semigroup object, groupoid object, etc.) in the category of quantum spaces. In this setting, the terms \textnormal{`abelian'} and \textnormal{`nonabelian'} will refer to the group structure of a quantum group space, and \textnormal{`commutative'} and \textnormal{`noncommutative'} will refer to its topology.
\end{defn}
With this definition, the main theme of the present article and its sequels is to simply replace $C^*$-algebras by the corresponding quantum spaces, and deduce results that can not be deduced, or even formulated, if we simply think of $C^*$-algberas as some abstract ``quantum spaces".  Following are some examples of such results. \newline

$\mathbf{(1)}$  Recall that the Gelfand duality for the unital  \emph{commutative} $C^*$-algebra generated by a normal operator $a$ leads to an integral fromula for the functional calculus of $a,$ and a special case of the formula is the Spectral Theorem \cite{dunford}. When $a$ is not assumed normal, our noncommutative Gelfand-Naimark duality  yields an integral formula (Theorem \ref{spectral}) for the noncommutative functional calculus for $a$, and a special case of this formula gives $a$ as an integral of an operator valued function with respect to a spectral measure on a space $Y$. This result is proved in Section \ref{functional}.

$\mathbf{(2)}$ In Section \ref{jordan-invariant}, we give an infinite analog of Jordan Canonical Form (Theorem \ref{jordan}) of an operator using which, we prove in Theorem \ref{invariant} the existence of nontrivial invaraint subspaces for an arbitrary bounded operator in a complex Hilbert space of dimension greater than one. The theorem has been hitherto proved in the case of normal operators, and several more general classes of operators (See \cite{kubrusly}, for example). The result is also true for compact operators \cite{kubrusly}, and holds trivially for arbitrary operators in \emph{nonseparable} Hilbert space. The case of arbitrary operator in infinite dimensional \emph{separable} Hilbert space had remained open. Our proof (Theorem \ref{invariant}) works for all operators. 

$\mathbf{(3)}$ In Section \ref{misc}, we describe two more applications, the proofs of which will appear in sequels of the present article: (i) The Pontryagin Duality theorem can be extended to abitrary locally compact groups (Theorem \ref{pontryagin}). The diagram in Section \ref{misc} gives a quick overview of this result. (ii) An extension of Stone's representation of Boolean algebras to orthomodular lattices. Strictly speaking, this is not an application of the results in the present article. Rather, it is an application of  the main idea of Theorem \ref{main} to an analogous problem in the field of Orthomodular lattices.\newline
\hbox{}\newline
Finally, the numbered remarks throughout the article point out how various results presented here reduce to standard results in commutative and/or finite dimensional cases.

Any phrase or symbol being defined will be typeset in bold face. Also, the use of the symbol $:=$ in an expression indicates that the left hand side is being defined.

\section{\normalsize{NONCOMMUTATIVE\ \ GELFAND-NAIMARK\ \ DUALITY}}\label{mainsection}

Let $A$ be a $C^*$-algebra. A state $\alpha$ on $A$ is pure if and only if the corresponding Gelfand-Naimark-Segal (GNS)  representation $\pi_{\alpha}$ is irreducible \cite{kadison}. 
\begin{defn}\label{equiv}
Note that $\overline{PS(A)},$ the weak*-closure of the set of pure states of $A,$ is compact Hausdorff,  and $0\in \overline{PS(A)}$ if and only if $A$ is without unit. Then  e define \[\bm{P(A)} := \overline{PS(A)}\setminus \{0\}.\] Note that $P(A)$ is a locally compact Hausdorff space, and is compact if and only if $A$ is unital. Now we say that $\alpha , \beta \in P(A)$ are  \textbf{equivalent} if the corresponding GNS representations are equivalent. 
We denote this equivalence relation by \[\bm{R(A)}\subset P(A)\times P(A)\]. 
\end{defn}
\begin{proposition} 
A $C^*$-algebra $A$ is commutative if and only if the equivalence relation $R(A)$ is discrete, i.e. all its equivalence classes are singleton sets. In this case, $R(A) = diag(P(A)\times P(A)).$
\end{proposition}
\textbf{Proof}:
The proof is trivial. At any rate, we are not going to use this result in what follows. Indeed it is an immediate corollary to Theorem \ref{main}. \ep

\begin{defn}\label{A-bar}
Let \[A \rightarrow C(P(A)): a\mapsto \overline{a}\] be the \textbf{functional representation} of $A$ given by \[\bm{\overline{a}(\alpha)} := \alpha(a).\] We denote the image of this representation by $\bm{\overline{A}}.$ 
\end{defn}
Now it is a wellknown fact that 
\begin{proposition}The following are equivalent:
\begin{enumerate}
\item  $A$ is commutative.
\item The map $a\mapsto \overline{a}$ is a $C^*$-algebra homomorphism.
\item The map $a\mapsto \overline{a}$ is onto, i.e. $\overline{A}= C(P(A)).$
\item $\overline{A}$ is a C*-subalgebra of $C(P(A)).$ 
\end{enumerate}
\ep
\end{proposition}

\begin{defn}\label{a-hat}
For each $a\in A,$ we define \[\bm{\widehat{a}}: \overline{A}\rightarrow \overline{A}\] by \[\bm{\widehat{a}(\overline{x})} := \overline{ax}\ \ \forall x\in A.\]
\end{defn}

\begin{remark}When $A$ is commutative, $\overline{A} =  C(P(A)),$ and $\widehat{a}(\overline{x}) = \overline{ax} = \overline{a}\ \overline{x},$ for all $x\in A.$ Thus $\widehat{a}: \overline{A}\rightarrow \overline{A}$ i.e. $\widehat{a}: C(P(A)) \rightarrow C(P(A))$ is the multiplication operator $\phi\mapsto \overline{a}\phi .$ This means that in the commutative case, $a \mapsto \widehat{a}\equiv \overline{a}$ is essentially the Gelfand transform.
\end{remark}

\begin{proposition}
The self-adjoint subspace $\overline{A}\subset C(P(A))$ generates $C(P(A))$ as a $C^*$-algebra.
\end{proposition}
\textbf{Proof:} Since $\overline{A}$ separates points of $P(A) \cup \{0\},$ so does the $C^*$-algbera  generated by $\overline{A}.$ Then the latter is equal to $C(P(A))$ by Stone-Weierstrass theorem.\ep 

Let $\bm{X}$ be a locally compact (Hausdorff) space 
Let $\bm{C(X)}$ be the $C^*$-algebra of continuous complex valued functions vanishing at infinity on $X$. Then the double dual of $C(X)$ is a von Neumann algebra, and its maximal ideal space $Y$ carries a canonical class of measures. What follows is independant of a choice of measure $m'$ in this class. Let $m$ be the image of $m'$ under the canonical onto map $Y\rightarrow X$. Then $m$ gives an embedding \[C(X)\hookrightarrow M(X): \phi\mapsto\phi\,dm,\] where $\bm{M(X)}$ is the Banach space of complex valued regular Borel measures on $X$. 

\begin{defn}\label{operator-measure}For a bounded linear operator \[\widehat{a}: C(X) \rightarrow M(X)\], define a regular complex valued Borel measure $\bm{d\widehat{a}}$ on $X \times X$  by the identity  \[\bm{\int{(f\otimes g)\hspace{1pt}d\widehat{a}}\,\,} :=  \int{g\ d (\widehat{a}f\!).}\] 
\end{defn}

\begin{defn}\label{contmeasure}
Let $C\subset C(X)$ be a self-adjoint closed subspace of $C(X)$ which separates points of $X$ and for each point in $x\in X$ there is a  $\phi\in C$ such that $\phi(x) \neq 0.$ Then $C$ generates $C(X)$ as a $C^*$algebra.
A measure $\mu \in M(X\times X)$ will be called a \textbf{C-measure} if $\mu = d(\widehat{a})$, for a linear operator $\widehat{a}: C \rightarrow C \subset C(X)\hookrightarrow M(X).$  For an equivalence relation $R(X)\subset X\times X$ let \[\bm{CO(R(X))}\] be the set of bounded operators on $C$ with support contained in $R(X)$, and \[\bm{CM(R(X))}\] be the set of $C$-measures with support contained in $R(X).$ \end{defn}

\begin{defn}\label{convolution}
Let $\mu , \nu \in CM(R(X))$ be given, and a Borel $E\subset R(X)$. Then define \[\bm{E'} := \{((x,y), (y,z))\in R(X)\times R(X): (x,z)\in E\}.\] Now define $\bm{\mu*\nu }$ by \[\bm{(\mu * \nu)}(E) := (\mu \times \nu )(E')\]. Then $\mu * \nu $, called the \textbf{convolution of $\mu$ and $\nu$},  is a measure contained in $CM(R(X))$.
\end{defn}

\begin{proposition}
The set $CM(R(X))$ is a $C^*$-algebra under convolution of measures. Furthermore, for $\widehat{a}, \widehat{b} \in CO(R(X)),$ let $\mu = d(\widehat{a}),$ and $\nu = d(\widehat{b})\in CM(R(X)),$ then $ d(\widehat{a}\widehat{b}) = \mu * \nu = d(\widehat{a}) * d(\widehat{b}),$ so that $CO(R(X)\rightarrow CM(R(X)$ defined by $\widehat{a}\mapsto d(\widehat{a})$ is an algebra isomorphism.\ep
\end{proposition}

Given $\widehat{a}\in CO(R(X)),$ let $d\widehat{a}\in CM(R(X))$ be the corresponding C-measure, and define $(\widehat{a})^*$ to be the element of $CO(R(X))$ corresponding to $(d\widehat{a})^*$.
Then \begin{proposition}\label{second}
 The map \[\widehat{a}\mapsto (\widehat{a})^*\] is an involution which makes $CO(R(X))$ a $C^*$-algebra, so that
$CO(R(X)) \cong CM(R(X))$ as $C^*$-algebras.\ep
\end{proposition}

Now let $X = P(A)$ defined \textit{supra},  and let $\bm{C} := \overline{A}\subset C(P(A))$ so that $\bm{CM(R(A))}$ is the convolution algebra of C-measures on $P(A)\times P(A)$ with support contained in $R(A)$, and $\bm{CO(R(A))}$ is the corresponding algebra of bounded operators $\overline{A} \rightarrow \overline{A}.$  

\begin{defn}For $a\in A$, let $\widehat{a}$ be the operator defined in Definition \ref{a-hat}, and let $d\widehat{a}$ be the corresponding measure as in Definition \ref{operator-measure}.
For every $a\in A$, we define \textbf{support of $\widehat{a}$}, \[\bm{Supp(\widehat{a})}:= Supp(d\widehat{a}),\] and define \[\bm{Supp(a)} := Supp(\widehat{a}).\]
\end{defn}

The following is the main result of this section:

\begin{theorem}[Noncommutative Gelfand-Naimark]\label{main}
For any $C^*$-algebra, $A$, the assignments $a\mapsto \widehat{a}\mapsto d\widehat{a}$ give the following $C^*$-isomorphisms:
\[A\hspace{4pt}\cong\hspace{4pt}CO(R(A))\hspace{4pt}\cong\hspace{4pt}CM(R(A)) .\]
\end{theorem}
\textbf{Proof:}
The second isomorphism is trivial (See Proposition \ref{second}). We prove the first isomorphism.
Let $\bm{\widehat{A}}$ be the image of the map $\widehat{\ \ }: A \rightarrow CO(R(A)): a \mapsto \widehat{a}.$
Clearly, $\widehat{\ \ }$ is linear, $\widehat{ab} = \widehat{a}\widehat{b},$ and $\widehat{a^*}$ = $(\widehat{a})^*.$

We now proceed to prove that the map $\widehat{\ \ }$  is (i) One-to-one, and (ii) Onto. Note that One-to-one implies  $A$ is isometric into $CO(R(A)):  A \rightarrow \widehat{A}\subset CO(R(A)).$

(i) Let $a\in A$.  Assume $\widehat{a} = 0.$ Then $\widehat{a}(\overline{x}) = 0, \forall x\in A.$ This implies that $\overline{ax} = 0, \forall x\in A.$ Consequently,  $\forall x\in A, $  \[ \forall \alpha\in P(A),\ \alpha(ax) = 0.\] This implies that $ax = 0, \ \forall x\in A,$ and hence $a = 0.$ Thus, for all \[a\in A, \ \ \widehat{a} = 0\  \Longrightarrow  \ a = 0.\] Thus, the map $\widehat{\ \ }$ is one to one, and hence an isometry, into $CO(R(A))$.

(ii)  Note that $R(A) =  R(CM(R(A)) = R(CM(R(\widehat{A})) = R(\widehat{A}).$ In particular, \[P(A) = P(\widehat{A}) = P(CM(R(A))) = P(CO(R(A)).\]  Now we know that $A$ separates points of $P(A)\cup \{0\}.$ Consequently, $\widehat{A}$ separates points of $P(A) \cup \{0\} = P(CO(R(A)))\cup\{0\}.$ Then, by the Noncommutative Stone-Weierstrass Theorem (See \cite{dixmier}, Corollary 11.5.2), $ \widehat{A}= CO(R(A)),$ and hence, $\widehat{ \ \ }:  A \rightarrow  CO(R(A))$ is a $C^*$-isomorphism.\ep

\begin{remark}For commmutative $A$,
$P(A)$ is the space of  pure states of $A$, i.e. the space of characters of $A.$
Then $R(A)$ is the diagonal of $P(A)\times P(A)$, and hence
$CO(R(A)) = CM(R(A)) =  C(P(A))$, so we recover the Gelfand-Naimark theorem (Theorem \ref{gelfand}).
\end{remark}
The theorem gives a quick proof of Dauns-Hoffman theorems \cite{dauns-hoffman}: (i) Representation of a $C^*$-algbera as continuous sections of a certain `sheaf' (ii) The center of a $C^*$-algebra $A$ is the algebra of continuous functions on the spectrum of $A$.

\section{\normalsize{FUNCTIONAL CALCULUS FOR HILBERT SPACE OPERATORS}}\label{functional}

In this section we present a noncommutative functional calculus for abitrary \linebreak bounded Hilbert space operators. 

Let  $\bm{a}$ be a bounded operator on a Hilbert space $\bm{H}$. Let $\bm{A}$ be the unital $C^*$-algebra generated by $\{1,a\}$,  and let $\bm{R(A)}$ be the equivalence relation defined by $A$ on $P(A)$.   When $a$ is normal, $A$ is commutative, and $R(A) = P(A) \cong \bm{\sigma(a)},$ the spectrum of $a$, and the spectral theory affords a formula for the functional calculus $L^{\infty}(\sigma(a))\rightarrow B(H): f\mapsto f(a)$   \[f(a) = \int_{\sigma(a)} f(z) \ dE,\] where $E$ is the spectral measure corresponding to $a$ \cite{dunford}. In particular,  for $f = z,$ the inclusion $z: \sigma(a)\!\hookrightarrow\!\mathbb{C},$ the formula reduces to \[a =  \int_{\sigma (a)} z \, dE ,\] which is, of course, the Spectral theorem.

In the general case, where $a$ is not assumed normal,  the functional calculus is noncommutative, and is given by a similar, albeit noncommutative, formula (Theorem \ref{spectral}). 

\begin{defn}
Let $\bm{q\!:\!P(A)\!\rightarrow\!Sp(A)},$ be the canonical quotient map corresponding to $R(A).$
Let $\mu \in CM(R(A)),$ and let $u$ be an $R(A)$-block, i.e. the cartesian product of an $R(A)$ equivalence class with itself. Then define the measure \[\bm{\mu_{u} \in CM(R(A))} \ \ \text{by}\ \ \bm{\mu_u(s)} := \mu(u\cap s),\] for all Borel subsets $s\subset R(A).$ Now for every $x\in Sp(A),$ let \[\bm{u_x} := q^{-1}(x) \times q^{-1}(x),\] an $R(A)$-block, and define an operator valued function \[\bm{\widehat{\mu}:  Sp(A)\rightarrow B(H)}\ \ \text{by}\ \ 
\bm{\widehat{\mu}(x)} := \mu_{u_x}(a).\] Let $\bm{\mu_z}$ be the measure in $CM(R(A))$ corresponding to the operator $a$, i.e. with $\mu_z(a) = a.$ The define the operator valued function \[\bm{\widehat{z}} := \widehat{\mu}_z.\] Let \[\bm{C(A)} := \{\widehat{\mu}: \mu\in CM(R(A))\}\]And let \[\bm{WM(R(A)) := CM(R(A))^{**}},\] be the double dual of $CM(R(A))$ considered as a subalgebra of $M(R(A)),$ and \[\bm{W(A)} := \{\widehat{\mu}: \mu\in WM(R(A))\}.\] Then $C(A) \subset W(A).$ 
If we define a norm on $W(A)$ by \[\bm{\|\widehat{\mu}\|} := sup\{\|\widehat{\mu}(x)\|: y\in Sp(A)\},\] then, it is easy to see that $\|\mu\| = \|\widehat{\mu}\|,$ so that \[W(A) \cong WM(R(A)), \ \ \text{and}\ \ C(A) \cong CM(R(A)).\]
Now for $\widehat{\mu}\in W(A),$ define \[\bm{\widehat{\mu}(a)} := \mu(a)\]. \end{defn}

\begin{theorem}[The functional calculus of a bounded operator]\label{spectral}
The functional calculus \[C(A)\rightarrow B(H)\] given by \[\widehat{\mu}\mapsto \widehat{\mu}(a)\] extends naturally to a functional calculus $W(A)\rightarrow B(H)$, and there exists a spectral measure on $Sp(A)$ such that $\widehat{\mu}\mapsto \widehat{\mu}(a)$ is given by the formula
\[\widehat{\mu}(a) = \int_{Sp(A)}\widehat{\mu} \  dE.\]
In particular,
\[ a = \int_{Sp(A)} \widehat{z}\ dE .\]
\end{theorem}
\textbf{Proof :}
For $\mu\in CM(R(A))$ and $g, h \in H$, define \[\bm{\alpha_{g,h}} := \ <\!\!\mu(a)g, h\!\!>\!\!.\] Then $\alpha_{g,h}$ is a bounded linear functional on $CM(R(A))$, with $ \|\alpha_{g,h}\| \leq \|g\| \ \|h\|$. Also, $(g, h) \mapsto \alpha_{g,h}$ is a sesquilinear form.  Extend $\alpha_{g,h}$ naturally to all $\mu\in WM(R(A)),$ by setting \[\bm{\alpha_{g,h}(\mu)} := \mu(\alpha_{g,h}).\] For any $\mu\in WM(R(A)),$ we can define 
$\bm{[\ , \ ]: H\times H \rightarrow \mathbb{C}}$ by \[\bm{[g, h]} := \alpha_{g,h}(\mu).\] Again, $[\ , \ ]$ is sesquilinear, and $|[g,h]| \leq \|\mu\| \ \|g\| \ \|h\|.$ Consequently, there is a unique $b\in B(H)$ such that $[g,h] = <\!\!b(g), h\!\!>,$ and $\|b\| \leq \|\mu\|.$ We now set \[\bm{\pi(\mu)} := b.\] Thus, $\|\pi(\mu)\| \leq \|\mu\|,$ and $\|\pi\| = 1$. \newline
\hbox{}\newline
Now we show that $\bm{\pi}: WM(R(A)) \rightarrow B(H)$ is a representation of $WM(R(A)).$

Let $\nu \in WM(R(A)) = CM(R(A))^{**}$, i.e.,  $f\mapsto \nu(f)$ is a bounded linear functional on $CM(R(A))^*$. Then $\{f\in CM(R(A)): \|f\| \leq \|\nu\| \}$ is dense  in the set $\{s\in WM(R(A)): \|s\| \leq \|\nu\|\},$ in the topology $\sigma(WM(R(A)), CM(R(A))^*\}.$ Hence there is a net $\mu_\alpha$ in $CM(R(A))$ such that $\|\mu_{\alpha}\| \leq \|\nu\|$, and $\forall f\in CM(R(A))^*$, $f(\mu_{\alpha}) \rightarrow f(\nu).$ Thus if  $\mu\in CM(R(A)),$ then $\pi(\nu *\mu) =  \pi(\lim \mu_\alpha * \mu) = \pi(lim (\mu_\alpha * \mu))  =  \lim(\pi(\mu_{\alpha})\pi(\mu)) = (\lim \pi(\mu_{\alpha})) \pi(\mu) = \pi(lim\ \mu_{\alpha})\ \pi(\mu) =   \pi(\nu)\ \pi(\mu).$ Now fix a $\mu\in WM(R(A)),$ and let $\mu_{\alpha}$ be a net in $CM(R(A))$ such that $\mu_{\alpha} \overset{W^*}{\rightarrow} \mu.$ Then $\pi(\nu*\mu) = lim\ \pi(\nu*\mu_{\alpha}) = lim \ (\pi(\nu)\ \pi(\mu_\alpha)) =  \pi(\nu) \ (lim \ \pi(\mu_{\alpha})) = \pi(\nu) \  \pi(\mu)$, where all limits are in $WOT$ topology.  Thus, $\forall \nu, \mu \in WM(R(A)),$  \[\pi(\nu * \mu) = \pi(\nu) \ \pi(\mu).\]

Now, $\pi(\nu)^* = \pi(\nu^*),$ for if $\{\nu_{\alpha}\}$ a net in $CM(R(A)),$ such that $\nu_\alpha \overset{W^*}{\rightarrow} \nu,$ then $\pi(\nu_{\alpha}) \overset{W^*}{\rightarrow} \pi(\nu),$ and hence $(\pi(\nu_\alpha))^*\overset{W^*}{\rightarrow} (\pi(\nu))^*.$ But  $(\pi(\nu_\alpha))^* = \pi(\nu_{\alpha}^*) \overset{W^*}{\rightarrow}\pi(\nu^*).$ Thus, \[\pi(\nu)^* = \pi(\nu^*).\]

Note that for $\mu \in CM(R(A)),$ we have $\pi(\mu) = \mu(a).$ Thus, $\pi: WM(R(A)) \rightarrow B(H)$ is an extension of $CM(R(A)) \rightarrow B(H): \mu \mapsto \mu(a).$ Then $\pi: W(A) \rightarrow B(H)$ defined by $\pi(\widehat{\mu}) := \pi(\mu)$ is a representation of $W(M).$

We now proceed to define the spectral measure $E$.
Let $1 \in CM(R(A))$ be the measure corresponding to the identity operator, and let $U$ be a Borel subset of $Y$. Now define  $1_U \in CM(R(A))$ by \[\bm{1_{U}} := \chi_{U}1.\]  Let $\Sigma$ be the Borel algebra of $Y$, define an operator valued measure $E$ on $\Sigma$, as follows. For all $U\in \Sigma,$
\[\bm{E(U)} := \pi(1_U),  \]
Then,
\begin{enumerate}
\item Since $1_U$ is a hermitian idempotent, $E(U)$ is a hermitian idempotent, i.e. $E(U)$ is a projection for all Borel $U\subset Y$.
\item $E(Y) = \pi(1) = 1$ and $E(\emptyset) = \pi(0) = 0$.

\item$E(U\cap V) = \pi(1_{U\cap V}) = \pi(1_U  1_V) = \pi(1_U)\ \pi(1_V) = E(U) E(V).$

\item Let $\{U_i\}_{1}^\infty$ be a sequence in $\Sigma$ such that $U_i \cap U_j = 0$ if $i \neq j$.  Then we show that $E(\bigcup_i(U_i)) =  \bigcup_i E(U_i).$   Now, $E$ is finitely additive by (ii) and (iii) above. Set $W_n = \bigcup_{n+1}^\infty U_i$. Then $\forall h\in H,$ we have $\|E(\bigcup_{i=1}^\infty U_i)h - \sum_{i=1}^{n} E(U_i)h \|^2 =\|E(\bigcup_{i=n+1}^\infty) h \|^2 = \|E(W_n) h\|^2  = \ <\!\!E(W_n)h,\ E(W_n)h\!\!>\ = \ <\!\!E(W_n)h, h\!\!> \ = \ <\!\!\pi(\mu_{W_n})h, \ h\!\!>\ = \alpha_{h,h}(\mu_{W_n}) =  \sum_{i=n+1}\alpha_{h,h}(U_n) \overset{WOT}{\rightarrow}0.$
\end{enumerate}
Thus, $E$ is a spectral measure.

Now we show that for each $\widehat{\mu}\in W(A)$ the formula $ \widehat{\mu}(a) = \int_{R(A)}\widehat{\mu}\ dE$ holds.
Fix a $\widehat{\mu}\in W(A).$ Let $\epsilon >  0$ be given. Then choose $U_i\in\Sigma, i = 1,...,n$ such that $U_i\cap U_j = 0$ if $i\neq j$, and $\bigcup_{i=1}^n U_i = Y $ and $sup\{\|\widehat{\mu}(y) - \widehat{\mu}(y')\|: y, y'\in U_i\}<  
 \epsilon$ for $1\leq i\leq n.$ Then, for any $y_i\in Y,$ $\|\widehat{\mu} - \sum_{i=1}^{n}\widehat{\mu }(y_i) \widehat{1}_{U_i}\| <   \epsilon. $ Now since $\|\pi\| = 1$, we have $\|\pi(\widehat{\mu}) - \sum_{i=1}^n \widehat{\mu}(y_i) \  \pi(\widehat{1}_{U_i}) \| \leq \|\pi\|  \ \|\widehat{\mu} - \sum_{i=1}^{n}\widehat{\mu }(y_i) \widehat{1}_{U_i}\| < \epsilon .$ 
Thus we have the formula \[\pi(\widehat{\mu}) = \int_{Y}\widehat{\mu} \ dE.\]\ep

\begin{remark}When $a\in B(H)$ is normal, $A$ is commutative,
$R(A) = P(A) =  \sigma (a) = Sp(A)$, so that $E$ and $\widehat{z}$ are the spectral measure and  the identity function respectively on $\sigma(a).$ COnsequently, the theorem reduces to the the functional calculus $L^{\infty}(\sigma(a))\rightarrow B(H),$ which includes the Spectral Theorem for normal operators.
\end{remark}

\section{\normalsize{INFINITE JORDAN CANONICAL FORM AND THE INVARIANT SUBSPACE THEOREM}}\label{jordan-invariant}
In this section we prove the Invariant Subspace Theorem (Theorem \ref{invariant}) by further refining the formula $a = \int_{Sp(A)} \widehat{z}\ dE$ of the preceding theorem. The following discussion leads up to this refinement which is, indeed, the infinite dimensional Jordan Canonical Form (Theoreom \ref{jordan}).

Let $\bm{\tau: P(A)\rightarrow \mathbb{C}}$ be the map given by $\bm{\tau(\alpha) := \alpha(a)}$, and  let $\bm{\Sigma(a)}$ be the image of $\tau.$ Then
\begin{proposition}\label{sigma}
\[\sigma(a) \subset \Sigma(a).\]\end{proposition}
\textbf{Proof:}
 Let $\lambda\in\sigma(a).$ Then $a-\lambda$ is not invertible. So $a-\lambda$ is not left invertible or not right invertible. If $a-\lambda$ is not left invertible, there exists a maximal left ideal $L$ of $A$ such that $a-\lambda \in L.$ Now there exists a pure state $\alpha\in P(A)$ such that $L=L_{\alpha}.$ Then $\alpha(a-\lambda) = 0,$ so $\alpha(a) = \lambda.$ On the other hand, if $a-\lambda$ is not right invertible, then there exists a maximal right ideal $R$ of $A$ such that $a-\lambda\in R.$ Then $L := R^*$ is a maximal left ideal, and $a^*-\overline{\lambda}\in L.$ Let $\alpha$ be a pure state such that $L=L_{\alpha}$. Then $\alpha(a^*- \overline{\lambda}) = 0.$ Now, since $\alpha$ is hermitian, $\overline{\alpha(a-\lambda)} = \alpha(a^*- \overline{\lambda}) = 0.$ So $\alpha(a) = \lambda .$
Thus $\lambda\in\Sigma(a).$\ep

Let $\bm{T(a)}$ be the equivalence relation defined on $P(A)$ by $\tau$. Let $\bm{R(a)}:= R(A) \cap T(a)$, the equivalence relation defined on $P(A)$. Let $r: P(A)\rightarrow Y$ be the quotient map corresponding to the relation $R(a)$. Then $\tau$ and $q$ factor through the quotient map $r$.  The following commutative diagram summarizes the situation.
\[
\xymatrix{
P(A)\ar@{->}[r]^{\tau}\ar@{->}[d]_{q}\ar@{->}[rd]^{r}&\Sigma(a)\ar@{<-}[d]^{z}&\sigma(a)\ar@{_(->}[l]\\
Sp(A) \ar@{<-}[r]_{b}& Y}
\]

For the rest of this section, it is extremely useful to keep this diagram in mind, and to think of the equivalence relations corresponding to the maps as consisting of blocks.

Note that corresponding to each $\lambda\in \Sigma(a)$, there may be  several $y \in z^{-1}(\lambda) \subset Y,$ which we will think of as different copies of $\lambda$, and for each of these there is an $R(a)$-class in $P(A).$  Thus, for each $\lambda\in \Sigma (a)$ there are several $R(a)$-blocks. This is roughly analogous to the fact that in the Jordan canonical form for operators in finite dimensional spaces, for a fixed $\lambda\in\sigma(a)$,  we may have several Jordan companion matrices $J_\lambda^k,$ of several different ranks $k,$ filling several disjoint diagonal square blocks. Thus, $z: Y \rightarrow \Sigma(a)$ can be thought of as a uniformization of $\Sigma (a)$.  Since,  by Proposition  \ref{sigma}, $\sigma(a)\subset\Sigma(a),$ the same considerations also apply to $\sigma(a)$. We will think of $R(a)$ as the scheme of blank blocks (i.e. $R(a)$-blocks) in our infinite dimensional case. Next step roughly amounts to filling these blocks with certain `Jordan matrices'. To this end, we recall the following standard definitions and facts:

\begin{defn} A bounded linear operator $n$ on a Hilbert space $H$ is called \textbf{quasi-nilpotent} if $\|n^i\|^{\frac{1}{i}}\rightarrow 0.$ Then $n$ is quasi-nilpotent if and only if $\sigma(a) = \{0\}.$ A bounded operator $s$ is said to be  \textbf{of scalar type} if it is of the form $\int_{\sigma(s)}\lambda\ dE$ where $E$ is the resolution of identity for $s$ \cite{dunford}.  A  \textbf{spectral operator} is, by definition \cite{dunford}, a bounded Hilbert space operator having a resolution of identity on the algebra of Borel subsets of its spectrum. Then, a bounded operator $b$ is spectral if and only it has a decomopsition $b = s + n,$ where $s$ is of scalar type, $n$ is quasi-nilpotent, and $sn=ns$ \cite{dunford}.  We make two more definitions: an operator $q$ will be called \textbf{quasi-scalar} if it is a scalar multiple of a projection. An operator $b$ will be called \textbf{strongly spectral} if $b = q+n$ where $q$ is quasi-scalar and $n$ is quasi-nilpotent. Note that a strongly spectral operator is spectral.
\end{defn}

Now the formula $a = \int_{Sp(A)} \widehat{z}\ dE,$ can be refined to yield a Jordan canonical form:
\begin{theorem}[Infinite Jordan Canonical Form]\label{jordan}
 Given a bounded operator $a$ on a Hilbert space $H,$
\[ a \ = \ \int_{Y} j \ dE ,\]
where $E$ is a spectral measure on $Y,$ $j: Y\rightarrow B(H)$ is an operator valued function such that $\forall y\in z^{-1}(\sigma(a)) \subset Y, \ j(y)$ is a strongly spectral operator.

Furthermore, there exists a spectral measure $E$ on $\Sigma(a),$ and an operator valued function $k:\Sigma(a)\rightarrow B(H),$ such that $\forall \lambda\in\sigma(a), k(\lambda)$ is a strongly spectral operator and \[a = \int_{\Sigma(a)} k \ dE.\]
\end{theorem}
\textbf{Proof:} Let $U$ be a Borel subset of $Y,$ and let $1_a\in CM(R(a))$ be the measure corresponding to the unit operator $1$. Then define the spectral measure $E$ on $Y$ exactly as in Theorem \ref{spectral}: $E(U) = \chi_{(U\times U) \cap R(a)}1_a,$ where $\chi_{(U\times U) \cap R(a)}$ is the characteristic function of $(U\times U) \cap R(a).$

Now we define the function $j: Y\rightarrow B(H)$ as follows.
For an $R(a)$-block $u$ define the measure $\bm{\mu_{u}}$ by \[\bm{\mu_{u}(s)} := \mu_z(u\cap s),\] for all Borel subsets $s\subset R(a).$ Now for every $y\in Y,$ define \[\bm{u_y} := r^{-1}(y)\times r^{-1}(y),\] an $R(a)$-block, and define an operator valued function $j:  Y\rightarrow B(H)$ by 
\[\bm{j(y)} := \mu_{u_y}(a).\]  Then the proof of the integral formula \[ a = \int_{Y} j \ dE \]is the same as the proof of Theorem \ref{spectral}.

 Let $\lambda\in\sigma(a)$, and $y\in Y$ such that $z(y) = \lambda.$  Then $\sigma(j(y)) = \{\lambda\},$ and $\sigma(j(y) - \lambda ) = \{0\}.$ Set $n(y) = j(y)-\lambda ,$ then $j(y) = \lambda  \,+\, n(y) = z(y) + n(y)$ with $\sigma(n(y))=\{0\},$ i.e. $n(y)$ is quasi-nilpotent. Thus $j(y)$ is a strongly spectral operator. 

Now we define a spectral measure $E$ on $\Sigma(a).$  
Let $V$ be a Borel subset of $\Sigma(a),$ and let $1_a\in CM(T(a))$ corresponding to the unit operator $1$. Also, define $V' = (\tau^{-1}(V)\times (\tau^{-1}(V)) \cap T(a).$ Then we define the spectral measure $E$ on $\Sigma$ by: \[\bm{E(V)} := \chi_{V'}1_a,\] where $\chi_{V'}$ is the characteristic function of $V'.$

Now we define the operator valued function $k:\Sigma(a)\rightarrow B(H).$ 

For a $T(a)$-block $u$, define the measure $\mu_{u}$ by \[\bm{u_{u}(s)} := \mu_z(u\cap s),\] for all Borel subsets $s\subset T(a).$ Now for every $\lambda\in \Sigma,$ define  \[\bm{u_{\lambda}} := \tau^{-1}(\lambda)\times \tau^{-1}(\lambda),\] a $T(a)$-block, and define an operator valued function $k:  \Sigma(a) \rightarrow B(H)$ by 
\[\bm{k(\lambda)} := \mu_{u_\lambda}(a).\]  

 Let $\lambda\in\sigma(a).$  Then $\sigma(k(\lambda)) = \{\lambda\},$ and $\sigma(k(\lambda) - \lambda ) = \{0\}.$ Set $n(\lambda) := k(\lambda)-\lambda ,$ then $k(\lambda) = \lambda  \,+\, n(\lambda),$ with $\sigma(n(\lambda))=\{0\},$ i.e. $n(\lambda)$ is quasi-nilpotent. Thus $k(\lambda)$ is a strongly spectral operator. 
Now, noting that $Supp(a)\subset T(a),$ the formula \[a = \int_{\Sigma(a)} k \ dE\]
holds in exactly the same way as the preceding formula:  $a = \int_{Y} j \ dE.$\ep

\begin{remark}
We remark that  for a \emph{finite} dimensional $H$  the formula $a = \int_{Y} j \ dE$ reduces to the Jordan Canonical Form. On the other hand, when $a$ is normal, the formula reduces to the Spectral Theorem.
\end{remark}
We will need the following well known fact for the proof of the Invariant Subspace Theorem ( Theorem \ref{invariant}).

\begin{lemma}\label{lemma}
Let $b$ and $c$ be nonzero operators on a Hilbert space $H$. If $bc = 0,$ then there exists a non-trivial subspace $S\subset H$ which is an invariant subspace for both $b$ and $c:$  i.e. $b(S)\subset S,$ and $c(S)\subset S.$
\end{lemma}
\textbf{Proof:}We include a proof of this fact for the sake of completeness. First some notation: for any $x\in B(H),$ we will denote by $\mathcal{R}(x)$ and $\mathcal{N}(x)$ the range and the kernel of $x$ respectively. Also, for $S\subset H,$ we will denote the closure of $S$ by $\overline{S}.$

Since we have assumed $b\neq 0,$ so $\mathcal{N}(b) \neq H.$ Now, since $bc = 0,$ we have $\mathcal{R}(c) \subset \mathcal{N}(b).$ Noting that $\mathcal{N}(b)$ is closed,  it follows that $H \neq \overline{\mathcal{R}(c)},$ the closure of $\mathcal{R}(c).$ Also, since $c \neq 0,$ we have $\mathcal{R}(c) \neq 0,$ and hence $\mathcal{N}(b) \neq 0.$ To summarize,
\[0 \neq \mathcal{N}(b) \neq H,\ \ \text{and}\ \ \ 0 \neq \overline{\mathcal{R}(c)} \neq H.\]
Now \[c(\mathcal{N}(b)) \subset c(H) = \mathcal{R}(c) \subset \mathcal{N}(b).\] 
Also, since $b(\mathcal{R}(c)) = 0,$ continuity of $b$ implies that \[b(\overline{\mathcal{R}(c)}) \subset\overline{b(\mathcal{R}(c))}\subset \overline{\mathcal{R}(c)}.\] Consequently, $\overline{\mathcal{R}(c)},$ and $\mathcal{N}(b)$ are nontrvial invariant subspaces for $b$ and $c$ respectively. On the other hand, for any $b$ and $c$, 
$\mathcal{N}(b)$ and $\overline{\mathcal{R}(c)}$ are invariant subspaces for $b$ and $c$ respectively.
Together, the last two assertion imply that $\mathcal{N}(b)$ and $\overline{\mathcal{R}(c)}$ are nontrivial invariant subspaces for \emph{both} $b$ and $c.$\ep

\begin{theorem}[Invariant Subspace Theorem]\label{invariant}
Every bounded operator on a complex Hilbert space $H$ with $dim(H)>1$ has a nontrivial invariant subspace.
\end{theorem}
\textbf{Proof:}
Let $a\in B(H).$  Then
\[a \ =\ \ \int_{\Sigma(a)} k\ dE \ \ = \ \ \int_{\sigma(a)} k\ dE \ \ \, + \ \ \int_{\Sigma(a)\setminus\sigma(a)} k\ dE.\]
Let $b = \int_{\sigma(a)} k\ dE$, and $c = \int_{\Sigma(a)\setminus\sigma(a)} k\ dE,$ so that $a = b \ +\ c$, and $bc = 0 = cb.$
Now we consider the following two cases.

(1) Assume $c=0$.
Then $a= b$ is a spectral operator, and as such has a nontrivial invariant subspace.  

(2) Assume $c \neq 0.$ Then by Lemma \ref{lemma}, $b$ and $c$ have a common nontrvial invariant subspace $S\subset H,$ so that $b(S)\subset S,$ and $c(S)\subset S.$ Now we consider the following two cases.\newline
(i) Assume $b(S) \neq 0.$ Then since $ab = b^2,$ we obtain \[a(b(S)) \subset b(b(S)) \subset b(S).\] Thus $b(S)$ is a nontrivial invariant subspace for $a.$\newline
(ii) Now assume that $b(S) = 0.$ Then \[a(S) = b(S) + c(S) = c(S) \subset S.\] Thus we see that $S$ is a nontrvial invariant subspace for $a$.

Together, (1) and (2) proves the theorem.\ep

\section{\normalsize{FURTHER\ \ APPLICATIONS}}\label{misc}
We describe here two more applications. Details, along with complete proofs, will appear in a sequel to the present article.

\subsection{\normalsize{Nonabelian Pontryagin Duality}}
Recall that the set $\widehat{G}$ of characters of a locally compact abelian group $G$ forms a locally compact abelian group and the celebrated Pontryagin duality theorem gives a natural isomorphism $G \cong \widehat{\widehat{G}}$.  We find that extending this theorem to nonabelian groups leads us to quantum group spaces as defined in Definition \ref{quantum space}: Given a locally compact group $G$ we define (see below) its dual to be a certain quantum group space $\widehat{G}$, which is a group if and only if $G$ is abelian. The classical dual of a possibly nonabelian $G,$ i.e. the set of equivalence classes of irreducible unitary representations of $G$,  is the quotient space corresponding to $\widehat{G}.$  In the abelian case, $\widehat{G}$ coincides with the classical dual.  This viewpoint inevitably leads to an extension of  the duality to \emph{quantum group spaces}. 

A quantum group space $G$ as defined in Definition \ref{quantum space}  is a group object in the category of quantum spaces. We attach a von Neumann bialgebra $K^*(G)$ to $G$.
 Now we can construct from the dual von Neumann bialgebra $\widehat{K^*(G)}$ a locally compact quantum space $\widehat{G}$, which has a multiplication  structure derived from the co-multiplication of $\widehat{K^*(G)}$.  This makes $\widehat{G}$ a quantum group space which we call \textbf{the dual quantum group space of $\bm{G}$}.  Following the same procedure, we construct a locally compact quantum group space $\widehat{\widehat{G}}$ from $\widehat{K^*(\widehat{G})}.$ Then we have the following generalization of classical Pontryagin duality .

\begin{theorem}[Pontryagin duality for quantum group spaces]\label{pontryagin}
For a quantum \linebreak group space $G$,   \[G \cong \widehat{\widehat{G}}.\] \ep
\end{theorem}
Recall from Definition \ref{quantum space} that the terms `abelian' and `nonabelian'  refer to the group structure of a quantum group space $G$, and `commutative' and `noncommutative' refer to the topology of $G$. Now,
let $G, H, K, N$ be quantum group spaces with the corresponding duals $\widehat{G}, \widehat{H}, \widehat{K}, \widehat{N}.$ Then the following table summarizes the various situations covered by Theorem  \ref{pontryagin}: 

\begin{center}\begin{tabular}[h]{|c|c|c|}
\hline \boldmath{$\downarrow$ Group-Space $\rightarrow$}& Commutative & Noncommutative\\ \hline  Abelian & \small{$G, \ \ \overset{}{\widehat{G}}$}&  
\small{$K,\ \ \overset{}{ \widehat{H}}$}\\\hline Nonabelian & \small{$H , \ \  \widehat{K}$}& 
\small{$N,  \ \ \overset{}{\widehat{N}}$}\\\hline
\end{tabular}
\end{center}
Thus, the dual $\widehat{G}$ of an abelian group $G$ is an abelian group;  for a nonabelian group $K$,  $\widehat{K}$ is an abelian noncommutative quantum group space, etc.  We note that the box containing $G, \widehat{G}$ is the classical Pontryagin duality. The boxes containing $H, \widehat{K}$ and  $K, \widehat{H}$ include nonabelian groups and abelian noncommutative group spaces, and finally the box containing $N, \widehat{N}$ cover nonabelian noncommutative quantum group spaces.

\subsection{\normalsize{Stone Duality for Noncommutative Boolean algebras, i.e. Orthomodular Lattices }} The ideas of Section  \ref{mainsection}  can be applied to Orthomodular Lattices (OML). In a certain precise sense, Boolean aglebras are commutative orthomodular lattices. Then, analogous to commutative $C^*$-Gelfand Duality, the Stone's representation theorem attaches a totally disconnected compact space to a complete Boolean algbera. Analogous to  Theorem \ref{main}, we can extend this theorem to general orthomodular lattices:

\begin{defn}\label{OML}
A set $L$ with operations $(\wedge , \vee , ', 0, 1)$  \textbf{is an orthomodular lattice} if  $\forall s, t \in L,$

\begin{enumerate}
\item $L(\wedge , \vee)$ is a lattice,
\item $(s')' = s ,$
\item $s \leqslant t  \Longrightarrow   t' \leqslant s' ,$ 
\item $ s \vee s' = 1,  s \wedge s'  =  0,$
\item$ s\leqslant t \Longrightarrow s \vee (s' \wedge t) = t$.
\end{enumerate}
The last condition is called the \textbf{orthomodularity} condition. 
\end{defn}
Lattices of projections in a $C^*$-algebras are the prime examples of OML.
Note that orthomodularity is a weakening of distributivity, so that distributive OML are simply Boolean algebras. For $a, b \in L,$ define  \[\bm{{a\Dot{\wedge}b}} := (a\vee b')\wedge b.\] Then \cite{beran} $L$ is a Boolean algebra if and only if $\forall \, a, b \in L,$\[a \ \Dot{\wedge}\  b = b\  \Dot{\wedge}\  a.\] In this sense, OML's are a noncommutative generalization of Boolean algebras. Also, if a $C^*$-algebra $A$ is generated by its lattice $L_A$ of projections (for example when $A$ is a von Neumann algebra), then $A$ is commutative if and only if the OML $L_A$ is commutative, i.e. a Boolean algebra. Now, elements of a Boolean algebra $B$ are represented by clopen subsets of a totally disconnected compact space---the maximal ideal space of $B$ (Stone's Theorem \cite{stone}).  As in the case of $C^*$-algebras, the geometric object corresponding to a (possibly noncommutative)\, OML is an equivalence relation on (or a quotient of) a totally disconnected compact space naturally associated with the lattice. Furthermore, an OML is Boolean if and only if this equivalence relation is discrete. In this case, one recovers Stone's theorem. The general case yields an OML analog of Dauns-Hoffman theorem---the Graves-Selesnick representation \cite{selesnick} of an OML as sections of a sheaf of (presumably simpler) OML's.

\end{document}